\documentclass[11pt,a4paper]{article}
\usepackage{color,amsmath,amssymb, lscape,graphicx,ifthen,pdfpages}
\usepackage[colorlinks=true,
linkcolor=webgreen,
filecolor=webbrown,
citecolor=webgreen]{hyperref}

\definecolor{webgreen}{rgb}{0,.5,0}\definecolor{webbrown}{rgb}{.6,0,0}

\newcommand{\seqnum}[1]{\href{http://oeis.org/#1}{\underline{#1}}}

\def\eg{{\it e.g.},\,}
\def\Eg{{\it E.g.},\,}
\def\ie{{\it i.e.},\,}
\def\Ie{{\it I.e.},\,}

\def\sspp{\,+\,}
\def\sspm{\,-\,}
\def\sspeq{\,=\,}
\def\sspdef{\, :=\,}
\def\sspneq{\,\neq\,}
\def\sspkl{\,<\,}
\def\sspgr{\,>\,}
\def\sspgeq{\,\geq\, }
 
\def\pn{\par\noindent}

\def\pbn{\par\bigskip\noindent}
\def\psn{\par\smallskip\noindent}

\def\sspequiv{\,\equiv\,}

\def\Beq{\begin{equation}}
\def\Eeq{\end{equation}}
\def\Beqarray{\begin{eqnarray}}
\def\Eeqarray{\end{eqnarray}}
\def\Eq#1{Eq.\,#1}

\def\sspin{\,\in\,}

\def\range#1#2{#1,\,\count15=#1 \advance\count15 by +1 \number\count15,\,...,\,#2\,}
\def\rangeinf#1{#1,\, \count16=#1 \advance\count16 by +1 \number\count16,\,...\,}

\normalbaselines
\begin{document}
\bibliographystyle{unsrt}
\rightline{Karlsruhe} \par\smallskip\noindent
\rightline{July 18}
\vbox {\vspace{6mm}}
\begin{center}
{\Large {\bf Cantor's List of Real Algebraic Numbers of Heights $\bf 1$ to $\bf 7$}}\\ [9mm]
Wolfdieter L a n g \footnote{ 
\href{mailto:wolfdieter.lang@partner.kit.edu}{\tt wolfdieter.lang@partner.kit.edu},\quad 
\url{http://www.itp.kit.edu/~wl}
                                          } \\[3mm]
\end{center}
\begin{abstract}
Cantor gave in his fundamental article \cite{Cantor} an elegant proof of the countability of real algebraic numbers based on a positive integer height, denoted by him as $N$, of integer and irreducible polynomials of given degree (denoted by him as $n$) with relative prime coefficients. The finite number of real algebraic numbers with given height he called $\varphi(N)$, and gave the first three instances.\pn
Here we give a systematic list for the real algebraic numbers of height, which we denote by $n$, for $n\sspeq 1,\, 2,\,...,\,7$ and polynomials of degree $k$.
\end{abstract}
\section {\bf Cantor's proof}
An {\sl algebraic number} is a number $\omega$ which satisfies a non-constant polynomial equation
\Beq
\sum_{j=0}^k\, a_j\, x^j \sspeq 0,\ {\text for}\  x\sspeq \omega.
\Eeq
where $k\sspin \mathbb N$, and the $a_j\sspin \mathbb Z$, for $j \sspin \{\range{0}{k}\}$. Without loss of generality one restricts to $a_k \sspgr 0$, and assumes that $\gcd(a_0,\,a_1,\,...,\,a_k)\sspeq 1$ (relative prime coefficients). In order to avoid multiple counting of numbers $\omega$ (in general complex, but in the following real) the polynomials have to be irreducible (they do not factorize over $\mathbb Z$). 
\psn
To prove that the set of such real algebraic numbers is countable, \ie that there is a bijective map between $\{\omega\}$ and $\mathbb N$, one has to list all such polynomials and restrict to their real roots, \ie discard all roots $\omega$ with non-vanishing imaginary part.  \psn
{\sl Cantor} \cite{Cantor} managed this by considering first polynomials with only non-negative coefficients $|a_j|$, and later allows all possible sign arrangements but keeping $a_k \sspgeq 1$.  Define
\Beq
K\sspdef \sum_{j=0}^{k-1}\,|a_j| \sspp a_k,
\Eeq 
hence $K$ is a positive integer. {\sl Cantor} introduced instead a {\sl height} $n\sspin \mathbb N$ satisfying, for $k\sspin \{\range{1}{n}\}$,
\Beq
n \sspeq K \sspp (k\sspm 1)\,.
\Eeq
 For each height $n$ there is a finite number of real roots (the polynomial of degree $k$ has at most $k$ distinct real roots), denoted $\Phi(n)$ (in \cite{Cantor} $\varphi(n)$) counted without roots which already appeared at lower heights. Thus one can give a list of the $\Phi(n,\,k)$ real algebraic numbers $\{\omega\}$ for each height $n$, and degree $k$ sorted in a specific way. Then these lists with possible $k$ values from $1$ to $n$ give, after concatenation, the list of $\Phi(n)$ entries for height $n$, for each $n$, for $n\sspgeq 1$. Thus each real algebraic number $\omega$, of degree $k$ determines the composition, \ie $K$, hence $n$. Therefore each $\omega$ has a unique address $c\sspin \mathbb N$ ($c$ for count), and for each $c\sspin N$ there is a unique $\omega$ by construction.
 \psn
Some remarks and implications depending on the height $n$ and degree $k$ follow.
\psn
{\bf (a)} To find the unsigned polynomials of height $n$, \ie the coefficients in \Eq{3}, one uses {\sl relatively prime compositions} (\ie partitions without regard of order and with only relatively prime parts). For the number of relatively prime partitions of $K$, with $K\sspin \mathbb N$, see \seqnum{A000837}$(K)$ (note that $\gcd(K)\sspeq K$). The corresponding number of relatively prime compositions is found in \seqnum{A000740}.\psn
These compositions (also the partitions) have positive integer parts, but here, except for $a_k$, the coefficients may also vanish. Therefore one uses sometimes also part $0$ in the later tables for the entries in the composition column. There the combinations will be recorded by falling powers of $x$ of the corresponding polynomial. As explained below such $0$s will only be shown if no confusion can arise.
\psn
{\bf (b)} Positive, negative or complex roots for signed polynomials.\pn
{\sl Descartes'} rule of signs (\url{https://mathworld.wolfram.com/DescartesSignRule.htm}) is useful to find the maximal number of (real) positive roots $\#r_+$  (counting multiplicities) or negative roots $\#r_-$ of a polynomial $P(k, x)$ (with real coefficients, but here with integer coefficients). Let the number of changes of signs in the coefficients be $S$, then  $\#r_+\sspin \{S,\, S-2,\, ..., 0\ {\rm or}\ 1\}$ The range of values of $\#r_-$ of $P(k,\,x)$ is obtained from the one of $\#r_+$ of $\pm P(k,\,-x)$.\pn
Concerning complex roots remember that for real polynomials they come in pairs of complex conjugate roots. \pn
\Eg For $P(3,\, x)\sspeq x^3 \sspp x^2\sspm x\sspm 1$ and for $-P(3,\, -x)\sspeq 
x^3 \sspm x^2\sspm x\sspp 1$ one finds $S\sspeq 1$ and $S\sspeq 2$, respectively. Hence there is one $r_+$ and either $\#r_-\sspeq 0$ (the case $r_+,\, c,\,\overline c$) or  $\#r_-\sspeq 2$, the case of three real roots. In fact, the given  $P(3,\,x)$ factorizes over the integers $P(3,\, x)\sspeq (x\sspm 1)\,(x\sspm 1)^2$, and the second case applies (with multiplicity counting) $1,\, -1,\, -1$. 
\psn
{\bf (c)} Case $k\sspeq 1$ which has $K\sspeq n$.
\pn
{\bf (i)} If $n\sspeq 1$ this gives the first real root $\omega(1)\sspeq 0$ because $a_0 \sspeq 0$, following from $a_1\sspgeq 1$. This leads to the unique polynomial $1\,x$, and the composition will be recorded as $[1]$, not $[1,0]$ (the $+$ sign of $a_k$ is never recorded).
\pn
{\bf (ii)} If $n\sspgeq 2$ one needs $|a_0|\sspgeq 1$ because if $a_0\sspeq 0$ the polynomial would give $n\,x$ (with no recorded $+$ sign, and $\gcd(n,\,0) \sspeq n\sspneq 1$), repeating the root $\omega(1)$ already found for height $n\sspeq 1$. This means that both parts $a_1$ and $|a_0|$ in the compositions are always non-vanishing, \ie the number of parts is $m\sspeq 2$. 
\pn
Example: For $n\sspeq 2$ one finds, with the composition $[1,1]$ with $K\sspeq 2$ and $m\sspeq 2$ the unsigned polynomial $x \sspp 1$ with root $-1$, and then the signed polynomial with root $+1$. These two roots will be recorded in this order: $\omega(2)\sspeq -1$, $\omega(3)\sspeq +1$.  
\psn
{\bf (d)} Case $k\sspeq 2$ with $K\sspeq n\sspm 1$. \ie $a_2\sspgeq 1$ and $|a_0|\sspgeq 1$ (no factorization).
\pn
{\bf (i)} Compositions with $a_1\sspeq 0$. For $n\sspeq 2$ this is impossible.
\pn 
Cases $n\sspgeq 3$: The candidates for relative prime compositions are ($a_1 = 0$ is omitted) $[n\sspm q,\, q\sspm 1]$ and $[q\sspm 1,n \sspm  q]$. For real roots and no factorization over the integers, with $q\sspgeq 2$ and $n\sspgeq q\sspp 1$, the unsigned version does not qualify because it leads to complex solutions.For the signed versions not both, $q\sspm 1$ and $n \sspm q$, can be squares (otherwise factorization occurs), and $\gcd(n \sspm q,\, q\sspm 1) \sspeq  \gcd(q\sspm 1,\, n \sspm q + (q\sspm 1)) \sspeq \gcd(q\sspm 1,n \sspm 1) \sspeq  1$. \Ie $\gcd(\overline{q\sspm 1}, n \sspm 1) \sspeq 1$, where $\overline m \sspeq$\seqnum{A007947}$(m)$ (the squarefree kernel of $m$). Note that always $n\sspm q \sspneq q\sspm 1$ (because equality is only possible for $q\sspeq 2$, and $n\sspeq 2\,q\sspm 1\sspeq 3$, but $2\sspm 1$ is a square, as well as $3\sspm 2$.
\pn 
Example $1$: $q\sspeq 2$, $\gcd(1,\,n-1) \sspeq 1$ for $n\sspgeq 3$, and $n\sspneq c^2\sspp 2$, for $c\sspgeq 1$. Thus for  $n\sspeq 3,\, 6,\, 11,\,18,\,...$ no such compositions qualify. $n\sspeq 4,\,5,\,7,\,8,\,9,
,10,\,12,\,...$\,.\pn
Example $2$: $q\sspeq 13$, $\overline{12}\sspeq 2\cdot3$,\, $n\sspgeq 14$: $\gcd(2\cdot3,\,n-1)\sspeq 1$, \ie $n$ is even and $n \sspequiv \{0,\,2\}\, ({\rm mod} \,3)$. Thus $n\sspeq {14,\, 18}\,({\rm mod}\, 6)$. The real roots come from signed compositions $[1,\, -12],\, [7,\,-12], ...$ and $[5,\, -12],\, [11,\, -12], ...$, and the ones with interchanged positive numbers.
\psn  
{\bf (ii)} Compositions with $|a_1|\sspgeq 1$.\pn
The compositions lead to the signed triples $[q,\,s1\,a_1,\,s2\,(n - 1 - (q + a_1))]$, with $q\sspeq a_2\sspgeq 1$, $a_1\sspgeq 1$, $n\sspgeq 2 + q + a_1$, and signs $s1$ and $s2$ from $\{+1,\,-1\}$.
\pn
Complex conjugate pairs of roots arise if $a_1^2 \sspm s2\,4\,q\,N(n,\,q,\,a_1) \sspkl 0$ with $N(n,\,q,\,a_1)\sspdef n - 1 - (q + a_1)\sspgeq 1$. \pn
If $a_1^2 \sspm s2\,4\,q\,N(n,\,q,\,a_1) \sspeq 0$ then factorization over $\mathbb Z$ of the corresponding polynomial appears.
\pn
If $a_1^2 \sspm s2\,4\,q\,N(n,\,q,\,a_1) \sspgr 0$ then also factorization over $\mathbb Z$ appears if $A(c)\sspdef (a_1^2\sspm c^2)/(4\,s2\,q)$ becomes an integer number, for integer $c\sspgeq  0$ and $A(c) \sspeq N(n,\,q,\,a_1)$. Otherwise the two real roots of the corresponding irreducible polynomial are the algebraic numbers
\Beq
x_{-}\sspeq -\frac{1}{2\,q}\,\left(s1\,a_1 \sspm \sqrt{a_1^2\sspm s2\,4\,q\,N(n,\,q,\,a_1)}\right),\ \rm{and}\  x_{+}\sspeq -\frac{1}{2\,q}\,\left (s1\,a_1 \sspp \sqrt{a_1^2\sspm s2\,4\,q\,N(n,\,q,\,a_1)}\right)\,.
\Eeq
Example: For $n\sspeq 4$ the composition is $[1,\,1,\,1]$. The signed triples with two real roots are $[1,\,+1,\,-1]$ and $[1,\,-1,\,-1]$. They are $-\varphi$ and $\varphi\sspm 1$, where $\varphi$ is the golden section \seqnum{A001622}, and $-(\varphi\sspm 1)$ and $\varphi$. repsectively.
This means that $\Phi(4,\,2)\sspeq 2\cdot2\sspp 2\cdot 2\sspeq 8$, including the above considered $a_1\sspeq 0$ cases. 
\psn
{\bf (e)} Case $k\sspeq n\sspgeq 1$ with $K\sspeq 1$.\pn
This can only appear for height $1$, when $a_0\sspeq 0$, giving the root $0\sspeq\omega(1)$. Thus $k$ runs at most from $1$ to $n-1$, for $n\sspgeq 2$.    
\psn
{\bf (f)} Case $k\sspeq n-1$, where $K\sspeq 2$, and this case cannot appear for $n\sspgeq 3$.\pn
This is because $a_n\sspgeq 1$ and $|a_0|\sspgeq 1$ (factorization otherwise), hence the polynomials are $x^{n-1}\sspp 1$ or $x^{n-1}\sspm 1$. The latter one is reducible ($n \sspgeq 3$), it factorizes into $(x\sspm 1)\,\sum_{j=0}^{n-2} {x^j}$. The first one has the negated roots of the second one if $n$ is even. If $n\sspeq 2\,q\sspp 1$, for $q\sspgeq 1$, then $\omega^q\sspeq \mp i$, hence $\omega$ is not real. Thus only degrees $k\sspin \{\range{1}{n-2}\}$ qualify if $n\sspgeq 3$.
\psn
{\bf (g)} Order of composition with $0$ parts (belonging to the same partition, not regarding their $0$ parts). \pn
These compositions are ordered such that the polynomials appear with falling powers of $x$.
\pn
The first two instances appears for $n\sspeq 5$ and $k\sspeq 3$ ($K\sspeq 3$), \ie $[1,\,1,\,0,\,1]$ and $[1,\,0,\,1,\,1]$. Both qualify for all $4$ sign assignments (the $+$ sign of $a_3$ is not recorded). \pn
Such $0$ parts are recorded whenever the number $m$ of (non-zero) parts of qualifying compositions of $K$ is less than $k+1$, except for $m\sspeq 1$ (only for $n\sspeq 1$ possible), and for all $m\sspeq 2$.
The number of (inner) $0$ parts is $z \sspeq k+1\sspm m$, for $m\sspgeq 3$. 
\psn
{\bf (h)} Order of the $\omega$s for each composition with given sign assignment (called signature).
The order of the $\omega$s increase. 
\psn
{\bf (i)} Order of compositions for given $n$ and $k$.\pn
The order for the list of the (qualifying) compositions of $K$, for $n \sspeq K\sspp k\sspm 1$, possibly with $0$ parts, obeys the order of the underlying relatively prime partitions of $K$ and given $n$ and $k$ (parts of partitions are ordered anti-lexicographically) with rising number of (positive) parts $m$. Compositions, not regarding $0$ parts, belonging to like $m$ partitions are ordered anti-lexicographically (\eg $[3,\,1,\,1]$ before $[2,\,2,\,1]$).  
\psn
{\bf (j)} Order of the signature for a given composition.
\pn
If there are different signatures for a given composition (possibly with zero parts) the order is in general with the first $\omega$s increasing. The exception is if it possible to arrange the order such that a $-|\omega|$ can be followed immediately by $|\omega|$. \Eg $\omega(40)$, then $\omega(41)$ (not $\omega(42)$).
\vfill\eject
\section {\bf The list of real algebraic numbers of height 1 to 7} 
The list of the first $291$ $\omega$s, numerated by $c\sspeq 1,\,, 2,\, ...,\,291$, is given in $7$ tables: 
\pn
{\bf Table 1}: $\omega(1)$,\, ...,\,$\omega(43)$\pn
{\bf Table 2}: $\omega(44)$,\, ...,\,$\omega(83)$\pn
{\bf Table 3}: $\omega(84)$,\, ...,\,$\omega(125)$\pn
{\bf Table 4}: $\omega(126)$,\, ...,\,$\omega(167)$\pn
{\bf Table 5}: $\omega(168)$,\, ...,\,$\omega(211)$\pn
{\bf Table 6}: $\omega(212)$,\, ...,\,$\omega(257)$\pn
{\bf Table 7}: $\omega(258)$,\, ...,\,$\omega(291)$.
\psn
The start is $\omega(0)\sspeq 0$ for $n\sspeq 1$ and $k\sspeq 1$, with $K\sspeq 1$, and the sign is not recorded (see remark {\bf (c) (i)}). \pn
We give some examples for the compositions illustrating the use (or omission) of $0$ parts, their order and the order of the signatures (signs) for a given composition.
\psn
{\bf Example 1:} Not recorded $0$ parts.\pn
For $n\sspeq 1$ see {\bf (c) (i)}.\pn
For all $n\sspgeq 4$ and $k\sspeq 2$ and partitions with number of parts $m\sspeq 2$, like $\omega(12)$.
\pn
{\bf Example 2:} Not recorded compositions.
\pn
(i)  For all $n\sspgeq2$ and $k\sspeq n$, because $K\sspeq 1$ (see remark $\bf (e)$).\pn
(ii) For all $\sspgeq 3$ and $k\sspeq n-1$. See remark $\bf (f)$. \pn
(iii) If for all signatures either the polynomials have pairs of complex conjugate roots or the polynomials factorize over the integers (the reducible case). \pn 
The instance of a composition where for all signatures the $\omega$s are complex cannot occur, because of {\sl Descartes'} rule of signs (see $\bf (b)$): start with all signs $+$ and assume only complex roots, then change the sign of the last coefficient, to obtain a positive real root.
\pn
For $n\sspeq 6$ and $k\sspeq 2$ the compositions $[4,\,1]$ and $[1,\,4]$ have two distinct pairs of  complex conjugate roots for sign $(+)$, and lead to factorization for sign $(-)$. \pn
For $n\sspeq 7$ and $k\sspeq 2$ the composition $[1,\,2,\,3]$ has a pair of complex conjugate roots for each of the signs $(+,+)$ and $(-,+)$, and factorization occurs for signs $(+,-)$  and $(-,-)$. The same applies to the compositions $[2,\,3,\,1]$, $[3,\,1,\,2]$ and $[3,\,2,\,1]$.
\pn
For $n\sspeq 7$ and $k\sspeq 4$ the compositions with two $0$ parts $[2,\,1,\,0,\,0,\,1]$, $[2,\,0,\,1,\,0,\,1]$ and $[2,\,0,\,0,\,1,\,1]$ lead each to a pair of complex conjugate $\omega$s for each of the signs $(+.+)$ and $(-,+)$, and factorization occurs for signs $(+,-)$ and $(-,-)$. The same applies to the compositions $[1,\,1,\,0,\,0,\,2]$, $[1,\,0,\,1,\,0,\,2]$ and $[1,\,0,\,0,\,1,\,2]$. 
\psn
{\bf Example 3:} Not recorded signatures for compositions.\pn
This happens if for recorded composition not all signatures are shown. Complex conjugat pairs may show up or factorization occurs.\pn
For many compositions of order $k\sspeq 2$ with two parts a complex conjugate pairs appear, like for sign $(+)$ with $n\sspeq 4$ compositions $[2,\,1]$ and $[1,\,2]$, or for $n\sspeq 5$ compositions $[3,\,1]$ and $[1,\,3]$.\pn
One complex conjugate pair may also appear for $k\sspeq 2$ compositions with more than two parts, like for $n\sspeq 5$ with composition $[1,\,2,\,1]$ and  signatures $(-,+)$ and $(+,+)$.\pn
For $k\sspeq 3$ one may find factorizations like for $n\sspeq 6$ and composition $[2,\,1,\,0,\,1]$ and signatures $(+,+)$ and $(-,-)$. For each case there is a real root found already for lower $n$ and, a complex conjugate pair.
\psn
{\bf Example 4:} Compositions with a given signature but not $k$ real solutions.\pn
Like for $n\sspeq 6$, $k\sspeq 4$, composition $[2,\,1]$ with sign $(-)$ and only two real roots (and a complex conjugate pair). 
\pbn
The number of $\omega$s for given height $n$ and degree $k$, \ie $\Phi(n,\, k)$, and the total number of real algebraic solution for $n$, \ie $\Phi(n)$, are given in {\sl Table 9}. The total numbers $\Phi(n)\sspeq$\seqnum{A362366}$(n)$ are, for $n\sspeq 1,\,2,...\,7$, $[1,\,2,\,4,\,12,\,28,\, 72,\, 172]$. \psn
\pbn

\pbn
\hrule
\pbn
OEIS\cite{OEIS} A-numbers:\,\seqnum{A000740},\,\seqnum{A000837},\,\seqnum{A007947},\,\seqnum{A364312},\,\seqnum{A364313},\,\seqnum{A364314},\,\seqnum{A364315},\,\seqnum{A364316}.\pn
For the other A-numbers see the seven tables.
\psn
Keywords: real algebraic numbers, height of algebraic numbers
\pn
MSC-Numbers: 11R04,\,13F20 
\pbn
\hrule
\pbn
\eject
\begin{center}
{\large {\bf Table 1: Real algebraic numbers $\bf \{{\boldsymbol\omega}(c)\}_{c=1}^{43}$ of height $\bf n$ and degree $\bf k$.}}
\par\smallskip\noindent
\begin{tabular}{|r|r|r|c|c|c|c|}\hline
&& && &&\\
$\bf c$ &  $\bf height\, n$ & $\bf degree\, k$ & \bf composition & $\bf signs $  & $\boldsymbol\omega$ & \bf A-number \\
&& && && \\ \hline\hline
$\bf 1$ &   $ 1 $  & $ 1 $ & $[1]$  &  $(\/)$ & {\color{blue}$\bf \ \ 0$} & $\seqnum{A001477}(0)$    \\
\hline
$\bf 2$ &   $ 2 $  & $ 1 $ & $[1,1]$  &  $(+)$ & {\color{blue}$\bf  -1$} & $- \seqnum{A001477}(1)$    \\
$\bf 3$ &   $   $  & $  $ &  $     $  &  $(-)$ & {\color{blue}$\bf +1$} & $+ \seqnum{A001477}(1)$    \\
\hline
$\bf 4$ &   $ 3 $  & $ 1 $ & $[2,1]$  &  $(+)$ & {\color{blue}$\bf -\frac{1}{2}$} & $-\seqnum{A020761}$    \\
$\bf 5$ &   $   $  & $   $ & $     $  &  $(-)$ & {\color{blue}$\bf +\frac{1}{2}$} & $+\seqnum{A020761}$    \\
$\bf 6$ &   $   $  & $   $ & $ [1,2]$ &  $(+)$ & {\color{blue}$\bf -2$} & $-\seqnum{A001477}(2)$    \\
$\bf 7$ &   $   $  &  $  $ & $     $  &  $(-)$ & {\color{blue}$\bf +2$} & $ +\seqnum{A001477}(2)$    \\
\hline
$\bf 8$ & $ 4 $  & $ 1 $  &  $[3,1]$  &  $(+)$ &  {\color{blue}$\bf -\frac{1}{3}$} & $-\seqnum{A010701}$  \\
$\bf 9$ & $   $  & $   $  &  $     $  &  $(-)$ &  {\color{blue}$\bf +\frac{1}{3}$} & $+\seqnum{A010701}$  \\
$\bf 10$& $  $  & $   $ &   $[1,3]$   &  $(+)$ &  {\color{blue}$\bf -3$} & $-\seqnum{A001477}(3)$  \\
$\bf 11$& $  $  & $   $ &  $    $   &   $(-)$  &  {\color{blue}$\bf +3$} & $+\seqnum{A001477}(3)$  \\
$\bf 12$ &   $  $  & $ 2  $ & $[2,1]$  &  $(-)$ & {\color{blue}$\bf -1/\sqrt{2}$} & $\mp \seqnum{A010503}$    \\
$\bf 13$ &   $  $  & $    $ & $     $  &  $   $ & {\color{blue}$\bf +1/\sqrt{2}$} & $\mp \seqnum{A010503}$    \\
$\bf 14$ &   $  $  & $   $  & $[1,2]$  &  $(-)$ & {\color{blue}$\bf -\sqrt{2}$} & $-\seqnum{A002193}$    \\
$\bf 15$ &   $  $  & $   $  & $     $  &  $    $ & {\color{blue}$\bf +\sqrt{2}$} & $+ \seqnum{A002193}$    \\
$\bf 16$ &   $  $  & $   $ & $[1,1,1]$  &  $(+,-)$ & {\color{blue}$\bf -\boldsymbol{\varphi} $} & $-\seqnum{A001622}$    \\
$\bf 17$ &  $  $  & $   $ & $        $  &  $     $ & {\color{blue}$\bf (\boldsymbol{\varphi}\sspm 1) $} & $\seqnum{A094214}$ \\
$\bf 18$ &  $  $  & $   $ & $        $  &  $(-,-)$ & {\color{blue}$\bf -({\boldsymbol\varphi}\sspm 1) $} & $-\seqnum{A094214}$ \\
$\bf 19$ &   $  $  & $   $ & $       $  &  $     $ & {\color{blue}$\bf {+\boldsymbol\varphi} $} & $\seqnum{A001622}$    \\
\hline
$\bf 20$ & $ 5 $  & $ 1 $ & $[4,1]$  &  $(+)$ &  {\color{blue}$\bf -\frac{1}{4}$} & $-\seqnum{A020773}$  \\
$\bf 21 $ & $  $  & $  $ & $    $    &  $(-)$ &  {\color{blue}$\bf +\frac{1}{4}$} & $+\seqnum{A020773}$  \\
$\bf 22$ & $ $  & $ $   & $[1,4]$    &  $(+)$ &  {\color{blue}$\bf -4$} & $- \seqnum{A001477}(4)$  \\
$\bf 23$ & $ $  & $ $   & $     $    &  $(-)$ &  {\color{blue}$\bf  +4$} & $ +\seqnum{A001477}(4)$  \\
$\bf 24$ & $ $  & $ $   & $[3,2]$    &  $(+)$ &  {\color{blue}$ -\bf\frac{2}{3}$} & $-\seqnum{A010722}$ \\
$\bf 25$ & $ $  & $ $   & $     $    &  $(-)$ &  {\color{blue}$ +\bf\frac{2}{3}$} & $+\seqnum{A010722}$ \\
$\bf 26$ & $ $  & $ $   & $[2,3]$    &  $(+)$ &  {\color{blue}$\bf -\frac{3}{2}$} & 
$-\seqnum{A152623}$  \\
$\bf 27$ & $ $  & $ $   & $     $  &  $(-)$ &  {\color{blue}$\bf +\frac{3}{2}$} & 
$ +\seqnum{A152623}$  \\
$\bf 28$ &   $  $  & $ 2  $ & $[3,1]$  &  $(-)$ & {\color{blue}$\bf -1/\sqrt{3}$} & $-\seqnum{A020760}$    \\
$\bf 29$ &   $  $  & $    $ & $     $  &  $   $ & {\color{blue}$\bf  +1/\sqrt{3}$} & $+\seqnum{A020760}$    \\
$\bf 30$ &   $  $  & $  $ & $[1,3]$  &    $(-)$ & {\color{blue}$\bf -\sqrt{3}$} & $-\seqnum{A002194}$  \\
$\bf 31$ &   $  $  & $  $ & $     $  &  $   $ & {\color{blue}$\bf  +\sqrt{3}$} & $+\seqnum{A002194}$  \\
$\bf 32$ &   $  $  & $  $ & $[1,2,1] $  &  $(+,-)$ & {\color{blue}$\bf -(\sqrt{2}\sspp 1)$} & $-\seqnum{A014176}$ \\
$\bf 33$ &   $  $  & $  $ & $        $  &  $     $ & {\color{blue}$\bf +\sqrt{2}\sspm 1$} & $+\seqnum{A188582}$  \\
$\bf 34$ &   $  $  & $  $ & $        $  &  $(-,-)$ & {\color{blue}$\bf -(\sqrt{2}\sspm 1)$} & $-\seqnum{A188582}$  \\
$\bf 35$ &   $  $  & $  $ & $        $  &  $   $ & {\color{blue}$\bf +(\sqrt{2}\sspp 1)$} & $+\seqnum{A014176}$ \\
$\bf 36$ &   $  $  & $ 3  $ & $[2,1]$  &  $(+)$ & {\color{blue}$\bf -(\frac{1}{2})^{\frac{1}{3}}$} & $-\seqnum{A270714}$    \\
$\bf 37$ &   $  $  & $   $ & $      $  &  $(-)$ & {\color{blue}$\bf +(\frac{1}{2})^{\frac{1}{3}}$} & $+\seqnum{A270714}$    \\
$\bf 38$ &   $  $  & $   $ & $[1,2]$  &  $(+)$ & {\color{blue}$\bf -2^{\frac{1}{3}}$} & $-\seqnum{A002580}$    \\
$\bf 39$ &   $  $  & $   $ & $     $  &  $(-)$ & {\color{blue}$\bf +2^{\frac{1}{3}}$} & $+\seqnum{A002580}$    \\
$\bf 40$ &   $  $  & $   $ & $[1,1,0,1]$  &  $(+,+)$ & {\color{blue}$\bf -1.4655712318...$} & $-\seqnum{A092526}$ \\
$\bf 41$ &   $  $  & $   $ & $      $  &   $(-,-)$ & {\color{blue}$\bf +1.4655712318...$} & $+\seqnum{A092526}$ \\
$\bf 42$ &   $  $  & $   $ & $      $  &  $(-,+)$ & {\color{blue}$\bf -0.7548776662...$} & $-\seqnum{A075778}$ \\
$\bf 43$ &   $  $  & $   $ & $      $  &   $(+,-)$ & {\color{blue}$\bf +0.7548776662...$} & $+\seqnum{A075778}$ \\
{\bf tbc} &  $\bf \cdots $  & $\bf \cdots $  &  $\bf \cdots $   &  $\bf \cdots $ &  $\bf \cdots $  & $\bf \cdots $ \\
\hline
\hline
\end{tabular}
\end{center}
\vfill
\eject
\begin{center}
{\large {\bf Table 2: Continued: Real algebraic numbers $\bf\{{\boldsymbol \omega}(c)\}_{c=44}^{83}$ of height $\bf n$ and degree $\bf k$.}}
\smallskip\noindent
\begin{tabular}{|r|r|r|c|c|c|c|}\hline
&& && && \\
$\bf c$ &  $\bf height\, n$ & $\bf degree\, k$ & \bf composition & $\bf signs $  & \bf $\boldsymbol{\omega}$ & \bf A-number \\
\hline\hline
$\bf 44$ &   $ 5$ cont.  & $  3 $ & $[1,0,1,1]$  &  $(-,+)$ & {\color{blue}$\bf -1.3247179572... $} & $-\seqnum{A060006}$    \\
$\bf 45$ &   $  $  & $   $ & $  $  &  $(-,-)$ & {\color{blue}$\bf +1.3247179572... $} & $+\seqnum{A060006}$    \\
$\bf 46$ &     & $   $ & $   $  &  $(+,+)$ & {\color{blue}$\bf -0.6823278038...$} & $-\seqnum{A263719}$ \\
$\bf 47$ &   $  $  & $   $ & $       $  &  $(+,-)$ & {\color{blue}$\bf +0.6823278038...$} & $+\seqnum{A263719}$ \\
\hline
$\bf 48$ &   $ 6 $  & $ 1 $ & $[5,1]$  &  $(+)$ & {\color{blue}$\bf  -\frac{1}{5}$} & $-\seqnum{A000038}$    \\
$\bf 49$ &   $   $  & $   $ & $     $  &  $(-)$ & {\color{blue}$\bf  +\frac{1}{5}$} & $+\seqnum{A000038}$    \\
$\bf 50$ &          &       & $[1,5]$  &  $(+)$ & {\color{blue}$\bf -5$} & 
$-\seqnum{A001477}(5)$    \\
$\bf 51$ &          &       & $     $  &  $(-)$ & {\color{blue}$\bf +5$} & 
$+\seqnum{A001477}(5)$    \\
$\bf 52$ &      & $ 2 $ & $[3,2]$  &  $(-)$ & {\color{blue}$\bf -\sqrt{\frac{2}{3}}$} & $-\seqnum{A157697}$    \\
$\bf 53$ &      & $   $ & $     $  &  $   $ & {\color{blue}$\bf +\sqrt{\frac{2}{3}}$} & $+\seqnum{A157697}$    \\
$\bf 54$ &      &       & $[2,3]$  &  $(-)$ & {\color{blue}$\bf -\sqrt{\frac{3}{2}} $} & $-\seqnum{A115754}$    \\
$\bf 55$ &      &       & $     $  &  $   $ & {\color{blue}$\bf +\sqrt{\frac{3}{2}} $} & $+\seqnum{A115754}$    \\
$\bf 56$ &      &       & $[3,1,1]$  &  $(+,-)$ & {\color{blue}$\bf -{\frac{1 + \sqrt{13}}{6}} $} & $-(\seqnum{A188943} - 1)$    \\
$\bf 57$ &      &       & $       $  &  $     $ & {\color{blue}$\bf +{\frac{-1 + \sqrt{13}}{6}} $} & $+\seqnum{A356033}$    \\
$\bf 58$ &      &       & $    $  &  $(-,-)$ & {\color{blue}$\bf -{\frac{  -1 + \sqrt{13}}{6}} $} & $-\seqnum{A356033}$    \\
$\bf 59$ &      &       & $       $  &  $     $ & {\color{blue}$\bf +{\frac{  1 + \sqrt{13}}{6}} $} & $+(\seqnum{A188943} - 1)$    \\
$\bf 60$ &      &       & $ [1,3,1]$  &  $(+,-)$ & {\color{blue}$\bf -{\frac{3 + \sqrt{13}}{2}}$} & $-\seqnum{A098316}$ \\
$\bf 61$ &      &       & $       $  &  $     $ & {\color{blue}$\bf {\frac{-3 + \sqrt{13}}{2}}$} & $+\seqnum{A085550}$ \\
$\bf 62$ &      &       & $   $  &  $(+,+)$ & {\color{blue}$\bf -(1 + {\boldsymbol \varphi})$} & $-\seqnum{A104457}$ \\
$\bf 63$ &      &       & $       $  &  $    $ & {\color{blue}$\bf -(2 - {\boldsymbol \varphi})$} & $-\seqnum{A132338}$ \\
$\bf 64$ &      &       & $       $  &  $(-,-)$ & {\color{blue}$\bf -{\frac{-3 + \sqrt{13}}{2}}$} &-$\seqnum{A085550}$ \\
$\bf 65$ &      &       & $       $  &  $     $ & {\color{blue}$\bf {\frac{3 + \sqrt{13}}{2}}$} & $+\seqnum{A098316}$ \\
$\bf 66$ &      &       & $       $  &  $(-,+)$ & {\color{blue}$\bf 2 -  {\boldsymbol\varphi}$} & +$\seqnum{A132338}$ \\
$\bf 67$ &      &       & $       $  &  $     $ &  {\color{blue}$\bf 1 + {\boldsymbol\varphi}$} & $+\seqnum{A104457}$ \\
$\bf 68$ &  $ $    &   $  $    & $[1,1,3]$  &  $(+,-)$ & {\color{blue}$\bf -{\frac{1 + \sqrt{13}}{2}}$} & $-\seqnum{A209927}$ \\
$\bf 69$ &      &               & $       $ &  $     $ & {\color{blue}$\bf +{\frac{-1 + \sqrt{13}}{2}}$} & $+\seqnum{A223139}$ \\
$\bf 70$ &  $ $    &   $ $      & $     $     &  $(-,-)$ & {\color{blue}$\bf -{\frac{-1 + \sqrt{13}}{2}}$} & $-\seqnum{A223139}$ \\
$\bf 71$ &         &          & $        $    &  $     $ & {\color{blue}$\bf +{\frac{1 + \sqrt{13}}{2}}$} & $+\seqnum{A209927}$ \\
$\bf 72$ &  $ $    &   $ $    & $ [2,2,1]$       &  $(+,-)$ & {\color{blue}$\bf -{\frac{1 + \sqrt{3}}{2}}$} & $-\seqnum{A332133}$ \\
$\bf 73$ &         &          & $     $       &     $     $ & {\color{blue}$\bf +{\frac{-1 + \sqrt{3}}{2}}$} & $+\seqnum{A152422}$ \\
$\bf 74$ &         &          & $     $       &  $(-,-)$ & {\color{blue}$\bf -{\frac{-1 + \sqrt{3}}{2}}$} & $-\seqnum{A152422}$ \\
$\bf 75$ &  $ $    &   $ $    & $ $           &  $     $ & {\color{blue}$\bf +{\frac{1 + \sqrt{3}}{2}}$} & $+\seqnum{A332133}$ \\
$\bf 76$ &  $ $    &   $ $    & $ [2,1,2]$    &  $(+,-)$ & {\color{blue}$\bf -{\frac{1 + \sqrt{17}}{4}}$} & $-\seqnum{A188934}$ \\
$\bf 77$ &         &          & $     $       &  $     $ & {\color{blue}$\bf +{\frac{-1 + \sqrt{17}}{4}}$} & $+(\seqnum{A188485}-1)$ \\
$\bf 78$ &  $ $    &   $ $    & $ $           &  $(-,-)$ & {\color{blue}$\bf-{\frac{-1 + \sqrt{17}}{4}}$} & $-(\seqnum{A188485}-1)$ \\
$\bf 79$ &         &          & $     $       &  $     $ & {\color{blue}$\bf +{\frac{1 + \sqrt{17}}{4}}$} & $+\seqnum{A188934}$ \\
$\bf 80$ &  $ $    &   $ $    & $[1,2,2] $    &  $(+,-)$ & {\color{blue}$\bf-(1 + \sqrt{3})$} & $-\seqnum{A090388}$ \\
$\bf 81$ &         &          & $     $       &  $     $ & {\color{blue}$\bf -1 + \sqrt{3}$} 
& $+\seqnum{A160390} $ \\
$\bf 82$ &         &          & $     $       &  $(-,-)$ & {\color{blue}$\bf 1 \sspm \sqrt{3}$} & $-\seqnum{A160390}$ \\
$\bf 83$ &  $ $    &   $ $    & $  $          &  $     $ & {\color{blue}$\bf 1 + \sqrt{3}$} & $+\seqnum{A090388}$ \\
\hline
{\bf tbc} &  $\bf \cdots $  & $\bf \cdots $  &  $\bf \cdots $   &  $\bf \cdots $ &  $\bf \cdots $  & $\bf \cdots $ \\
\hline
\end{tabular}
\end{center}
\vfill
\eject
\begin{center}
{\large {\bf Table 3: Continued: Real algebraic numbers $\bf \{{\boldsymbol \omega}(c)\}_{c=84}^{125}$  of height $\bf n$ and degree $\bf k$.} }
\par\smallskip\noindent
\begin{tabular}{|r|r|r|c|c|c|c|}\hline
&& && &&\\
$\bf c$ &  $\bf height\, n$ & $\bf degree\, k$ & \bf Composition & $\bf signs $  & $\boldsymbol \omega$  & \bf A-number \\
&& && && \\ \hline\hline
$\bf 84$ &   $ 6$ cont.  & $ 3 $ & $[3,1]$          &  $(+)$ & {\color{blue}$\bf  -{(\frac{1}{3})}^{\frac{1}{3}}$} & $-\seqnum{A072365}$    \\
$\bf 85$ &   $   $  & $   $ & $     $          &  $(-)$ & {\color{blue}$\bf +{(\frac{1}{3})}^{\frac{1}{3}}$} & $+\seqnum{A072365}$    \\
$\bf 86$ &          &       & $[1,3]$          &  $(+)$ & {\color{blue}$\bf -{3}^{\frac{1}{3}}$}  & $-\seqnum{A002581}$    \\
$\bf 87$ &          &       & $     $          &  $(-)$ & {\color{blue}$\bf +{3}^{\frac{1}{3}}$}  & $+\seqnum{A002581}$    \\
$\bf 88$ &  $ $    &   $ $    & $ [2,1,0,1]$   &  $(-,+)$ & {\color{blue}$\bf -0.6572981061...$} & $-\seqnum{A089826}$ \\
$\bf 89$ &  $ $    &   $ $    & $   $          &  $(+,-)$ & {\color{blue}$\bf +0.6572981061...$} & $+\seqnum{A089826}$ \\
$\bf 90$ &  $ $    &   $ $    & $ [2,0,1,1]$   &  $(+,+)$ & {\color{blue}$\bf -0.5897545123...$} & $-\seqnum{A356031}$ \\
$\bf 91$ &  $ $    &   $ $    & $   $          &  $(+,-)$ & {\color{blue}$\bf +0.5897545123...$} & $+\seqnum{A356031}$ \\
$\bf 92$ &  $ $    &   $ $    & $ [1,2,0,1]$   &  $(+,+)$ & {\color{blue}$\bf -2.2055694304...$} & $-\seqnum{A356035}$ \\
$\bf 93$ &  $ $    &   $ $    & $   $          &  $(-,-)$ & {\color{blue}$\bf +2.2055694304...$} & $+\seqnum{A356035}$ \\
$\bf 94$ &  $ $    &   $ $    & $ [1,0,2,1]$   &  $(+,+)$ & {\color{blue}$\bf -0.4533976515...$} & $-\seqnum{A272874}$ \\
$\bf 95$ &  $ $    &   $ $    & $   $          &  $(-,-)$ & {\color{blue}$\bf +0.4533976515...$} & $+\seqnum{A272874}$ \\
$\bf 96$ &  $ $    &   $ $    & $ [1,1,0,2]$   &  $(+,+)$ & {\color{blue}$\bf -1.6956207695...$} & $-\seqnum{A289265}$ \\
$\bf 97$ &  $ $    &   $ $    & $   $          &  $(-,-)$ & {\color{blue}$\bf +1.6956207695...$} & $+\seqnum{A289265}$ \\
$\bf 98$ &  $ $    &   $ $    & $ [1,0,1,2]$   &  $(-,+)$ & {\color{blue}$\bf -1.5213797068...$} & $-\seqnum{A356030}$ \\
$\bf 99$ &  $ $    &   $ $    & $   $          &  $(-,-)$ & {\color{blue}$\bf +1.5213797068...$} & $+\seqnum{A356030}$ \\
$\bf 100$ &  $ $    &   $ $    & $ [1,1,1,1]$   &  $(+,-,+)$ & {\color{blue}$\bf -1.8392867552...$} & $-\seqnum{A058265}$ \\
$\bf 101$ &  $ $    &   $ $    & $   $          &  $(-,-,-)$ & {\color{blue}$\bf +1.8392867552...$} & $+\seqnum{A058265}$ \\
$\bf 102$ &  $ $    &   $ $    & $   $          &  $(-,+,+)$ & {\color{blue}$\bf -0.5436890126...$} & $-\seqnum{A192918}$ \\
$\bf 103$ &  $ $    &   $ $    & $   $          &  $(+,+,-)$ & {\color{blue}$\bf +0.5436890126...$} & $+\seqnum{A192918}$ \\
\hline
$\bf 104$ &   $ 6  $  & $ 4  $ & $[2,1]$  &  $(-)$ & {\color{blue}$\bf -(\frac{1}{2})^{\frac{1}{4}}$} & $-\seqnum{A228497}$    \\
$\bf 105$ &   $    $  & $    $ & $     $  &  $   $ & {\color{blue}$\bf +(\frac{1}{2})^{\frac{1}{4}}$} & $+\seqnum{A228497}$    \\
$\bf 106$ &  $ $  & $  $    & $[1,2]$     &  $ (-)  $ & {\color{blue}$\bf -2^{\frac{1}{4}}$} & $ -\seqnum{A010767}$    \\
$\bf 107$ &  $ $  & $  $    & $     $     &  $  $ & {\color{blue}$\bf +2^{\frac{1}{4}}$} & $+\seqnum{A010767}$    \\
$\bf 108$ &  $ $  & $  $    & $[1,1,0,0,1]$   &  $(+,-)$ & {\color{blue}$\bf -1.3802775690...$} &  $-\seqnum{A086106}$ \\
$\bf 109$ &  $ $  & $  $    & $           $   &  $     $ & {\color{blue}$\bf +0.8191725133...$}    &  $+\seqnum{A230151}$ \\
$\bf 110$ &  $ $  & $  $    & $       $   &  $(-,-)$ & {\color{blue}$\bf -0.8191725133...$} &  $-\seqnum{A230151}$ \\
$\bf 111$ &  $ $  & $  $    & $       $   &  $     $ & {\color{blue}$\bf +1.3802775690...$}    &  $+\seqnum{A086106}$ \\
$\bf 112$ &  $ $  & $  $    & $ [1,0,1,0,1]$   &  $(-,-)$ & {\color{blue}$\bf -\sqrt{{\boldsymbol\varphi}} $} &  $-\seqnum{A139339}$ \\
$\bf 113$ &  $ $  & $  $    & $       $   &  $     $ & {\color{blue}$\bf \sqrt{{\boldsymbol\varphi}}$}    &  $+\seqnum{A139339}$ \\
$\bf 114$ &  $ $  & $  $    & $            $   &  $(+,-)$ & {\color{blue}$\bf -\sqrt{{\boldsymbol\varphi} \sspm 1}$} &  $-\seqnum{A197762}$ \\
$\bf 115$ &  $ $  & $  $    & $           $   &  $     $ & {\color{blue}$\bf \sqrt{{\boldsymbol\varphi} \sspm 1} $}    &  $+\seqnum{A197762}$ \\
$\bf 116$ &  $   $  & $  $   & $[1,0,0,1,1]$   &  $(+,-)$ & {\color{blue}$\bf -1.2207440846...$} &  $-\seqnum{A060007}$ \\
$\bf 117$ &  $ $  & $  $    & $           $   &  $     $ & {\color{blue}$\bf +0.7244919590... $}    &  $+\seqnum{A356032}$ \\
$\bf 118$ &  $ $  & $  $    & $       $   &  $(-,-)$ & {\color{blue}$\bf -0.7244919590... $} &  $-\seqnum{A356032}$ \\
$\bf 119$ &  $ $  & $  $    & $       $   &  $     $ & {\color{blue}$\bf +1.2207440846...$}    &  $+\seqnum{A060007}$ \\
\hline
$\bf 120$ & $ 7 $  & $ 1 $ & $[6,1]$  &  $(+)$ & {\color{blue}$\bf -{\frac{1}{6}}$} & $-\seqnum{A020793}$   \\
$\bf 121$ & $   $  & $   $ & $     $  &  $(-)$ & {\color{blue}$\bf +{\frac{1}{6}}$} & $+\seqnum{A020793}$   \\
$\bf 122$ &   $  $  & $  $ & $[1,6]$  &  $(+)$ & {\color{blue}$\bf -6$} & $-\seqnum{A001477}(6)$ \\
$\bf 123$ &   $  $  & $  $ & $     $  &  $(-)$ & {\color{blue}$\bf +6$} & $+\seqnum{A001477}(6)$ \\
$\bf 124$ & $  $  & $  $ & $[5,2]$  &  $(+)$ &  {\color{blue}$\bf -{\frac{2}{5}}$} & $-(\seqnum{A274981}-1)$  \\
$\bf 125$ & $  $  & $  $ & $     $  &  $(-)$ &  {\color{blue}$\bf +{\frac{2}{5}}$} & $+(\seqnum{A274981}-1)$  \\
{\bf tbc} &  $\bf \cdots $  & $\bf \cdots $  &  $\bf \cdots $   &  $\bf \cdots $ &  $\bf \cdots $  & $\bf \cdots $ \\
\hline
\end{tabular}
\end{center}
\vfill
\eject
\begin{center}
{\large {\bf Table 4: Continued: Real algebraic numbers $\bf\{{\boldsymbol \omega}(c)\}_{c=126}^{167}$ of height $\bf n$ and degree $\bf k$.}  }
\psn
\begin{tabular}{|r|r|r|c|c|c|c|}\hline
$\bf c$ &  $\bf height\, n$ & $\bf degree\, k$ & \bf Composition & $\bf signs $  & $\boldsymbol \omega$  & \bf A-number 
\\ \hline\hline
$\bf 126$ &   $ 7$ cont.  & $ 1$ cont. & $[2,5]$  &  $(+)$ & {\color{blue}$\bf -{\frac{5}{2}}$} & $-10*\seqnum{A020773}$    \\
$\bf 127$ &   $  $  & $  $ & $     $  &  $(-)$ & {\color{blue}$\bf +{\frac{5}{2}}$} & $+10*\seqnum{A020773}$    \\
$\bf 128$ &   $  $  & $  $ & $[4,3]$  &  $(+)$ & {\color{blue}$\bf -{\frac{3}{4}}$} & $-\seqnum{A152627}$    \\
$\bf 129$ &   $  $  & $  $ & $     $  &  $(-)$ & {\color{blue}$\bf +{\frac{3}{4}}$} & $+\seqnum{A152627}$    \\
$\bf 130$ &   $  $  & $  $ & $[3,4]$  &  $(+)$ & {\color{blue}$\bf -{\frac{4}{3}}$} & $-\seqnum{A122553}$    \\
$\bf 131$ &   $  $  & $  $ & $     $  &  $(-)$ & {\color{blue}$\bf +{\frac{4}{3}}$} & $+\seqnum{A122553}$    \\
\hline
$\bf 132$ & $  $  & $ 2 $  &  $[5,1]$  &  $(+)$ &  {\color{blue}$\bf  -(-1 \sspp 2\,{\boldsymbol\varphi})$} & $-\seqnum{A002163}$  \\
$\bf 133$ & $  $  & $   $  &  $    $  &  $(-)$ &  {\color{blue}$\bf +(-1 \sspp 2\,{\boldsymbol\varphi})$} & $+\seqnum{A002163}$  \\
$\bf 134$ & $  $  & $   $  & $[1,5]$  &  $(+)$ & {\color{blue}$\bf  -(-1 \sspp 2\,{\boldsymbol\varphi})/5 $} & $-\seqnum{A002163}/5$    \\
$\bf 135$ & $  $  & $   $  & $     $  &  $(-)$ & {\color{blue}$\bf  +(-1 \sspp 2\,{\boldsymbol\varphi})/5 $} & $+\seqnum{A002163}/5$    \\
$\bf 136$ &  $$    &   $  $    & $[4,1,1]$  &  $(+,-)$ & {\color{blue}$\bf -{\frac{1 + \sqrt{17}}{8}}$} & $-(\seqnum{A189038}\sspm 1)$ \\
$\bf 137$ &  $$    &   $  $    & $       $  &  $     $ & {\color{blue}$\bf {\frac{-1 \sspp \sqrt{17}}{8}}$} & $+\seqnum{A358945}$ \\
$\bf 138$ &  $$    &   $  $    & $       $  &  $(-,-) $ & {\color{blue}$\bf -{\frac{-1 \sspp \sqrt{17}}{8}}$} & $-\seqnum{A358945}$ \\
$\bf 139$ &  $$    &   $  $    & $       $  &  $   $ & {\color{blue}$\bf +{\frac{1 + \sqrt{17}}{8}}$} & $+\seqnum{A189038}\sspm 1$ \\
$\bf 140$ &  $$    &   $  $    & $ [1,4,1]$     &  $(+,-)$ &\bf {\color{blue}$\bf -(1\sspp 2\,{\boldsymbol\varphi})$} & $-\seqnum{A098317}$ \\
$\bf 141$ &  $$    &   $  $    & $       $  &  $     $ & {\color{blue}$\bf +(-3 \sspp 2\,{\boldsymbol\varphi})$} & $+(\seqnum{A134972}\sspm 1)$ \\
$\bf 142$ &  $$    &   $  $    & $        $  &  $(+,+)$ & {\color{blue}$\bf -(2 \sspp \sqrt{3})$} & $-\seqnum{A019973}$ \\
$\bf 143$ &  $$    &   $  $    & $       $  &  $     $ & {\color{blue}$\bf -(2 \sspm \sqrt{3})$} & $-\seqnum{A019913}$ \\
$\bf 144$ &  $$    &   $  $    & $       $  &  $(-,-)$ & {\color{blue}$\bf -(-3 \sspp 2\,{\boldsymbol\varphi})$} & $ -(\seqnum{A134972}\sspm 1)$ \\
$\bf 145$ &  $$    &   $  $    & $    $     &  $     $ & {\color{blue}$\bf + (1\sspp 2\,{\boldsymbol\varphi} )$} & $+\seqnum{A098317}$ \\
$\bf 146$ &  $$    &   $  $    & $       $  &  $(-,+)$ & {\color{blue}$\bf +(2 \sspm \sqrt{3})$} & $+\seqnum{A019913}$ \\
$\bf 147$ &  $$    &   $  $    & $    $     &  $     $ & {\color{blue}$\bf +(2 + \sqrt{3})$} & $-\seqnum{A019973}$ \\
$\bf 148$ &  $$    &   $  $    & $[1,1,4]$  &  $(+,-)$ & {\color{blue}$\bf -{\frac{1 \sspp \sqrt{17}}{2}}$} & $-\seqnum{A222132}$ \\
$\bf 149$ &  $$    &   $  $    & $       $  &  $     $ & {\color{blue}$\bf +{\frac{-1 \sspp \sqrt{17}}{2}}$}   & $+\seqnum{A222133 }$ \\
$\bf 150$ &  $$    &   $  $    & $       $  &  $(-,-) $ & {\color{blue}$\bf -{\frac{-1 \sspp \sqrt{17}}{2}}$} & $-\seqnum{A222133}$ \\
$\bf 151$ &  $$    &   $  $    & $       $  &  $   $ & {\color{blue}$\bf +{\frac{1 \sspp \sqrt{17}}{2}}$} & $+\seqnum{A222132}$ \\
$\bf 152$ &  $$    &   $  $    & $[2,3,1]$  &  $(+,-)$ & {\color{blue}$\bf -{\frac{3 \sspp \sqrt{17}}{4}}$} & $-\seqnum{A188485}$ \\
$\bf 153$ &  $$    &   $  $    & $       $  &  $     $ & {\color{blue}$\bf +{\frac{-3 \sspp \sqrt{17}}{4}}$}   & $+\seqnum{A246725}$ \\
$\bf 154$ &  $$    &   $  $    & $       $  &  $(-,-) $ & {\color{blue}$\bf -{\frac{-3 \sspp \sqrt{17}}{4}}$} & $-\seqnum{A246725}$ \\
$\bf 155$ &  $$    &   $  $    & $       $  &  $   $ & {\color{blue}$\bf +{\frac{3 \sspp \sqrt{17}}{4}}$} & $+\seqnum{A188485}$ \\
$\bf 156$ &  $$    &   $  $    & $[1,3,2]$  &  $(+,-)$ & {\color{blue}$\bf -{\frac{3 \sspp \sqrt{17}}{2}}$} & $-\seqnum{A178255}$ \\
$\bf 157$ &  $$    &   $  $    & $       $  &  $     $ & {\color{blue}$\bf +{\frac{-3 \sspp \sqrt{17}}{2}}$}   & $+(\seqnum{A178255}\sspm 3)$ \\
$\bf 158$ &  $$    &   $  $    & $       $  &  $(-,-) $ & {\color{blue}$\bf -{\frac{-3 \sspp \sqrt{17}}{2}}$} & $-(\seqnum{A178255}\sspm 3)$ \\
$\bf 159$ &  $$    &   $  $    & $       $  &  $   $ & {\color{blue}$\bf +{\frac{3 \sspp \sqrt{17}}{2}}$} & $+\seqnum{A178255}$ \\
\hline
$\bf 160$ & $  $  & $ 3 $  &  $[4,1]$  &  $(+)$ &  {\color{blue}$\bf -\frac{2^{1/3}}{2} $ } & $-\seqnum{A235362}$  \\
$\bf 161$ & $  $  & $   $  &  $     $  &  $(-)$ &  {\color{blue}$\bf +\frac{2^{1/3}}{2} $ } & $+\seqnum{A235362}$  \\
$\bf 162$ & $  $  & $   $  &  $[1,4]$  &  $(+)$ & {\color{blue}$\bf -\frac{2}{2^{1/3}}$} & $-\seqnum{A005480}$    \\
$\bf 163$ & $  $  & $   $  &  $    $  &  $(-)$ & {\color{blue}$\bf  +\frac{2}{2^{1/3}} $} & $+\seqnum{A005480}$    \\
$\bf 164$ & $  $  & $  $  &  $[3,2]$  &  $(+)$ &  {\color{blue}$\bf -\frac{(18)^{1/3}}{3} $ } & $-\seqnum{A358943}$  \\
$\bf 165$ & $  $  & $  $  &  $     $  &  $(-)$ &  {\color{blue}$\bf +\frac{(18)^{1/3}}{3} $ } & $+\seqnum{A358943}$  \\
$\bf 166$ & $  $  & $   $  &  $[2,3]  $    &  $(+)$ & {\color{blue}$\bf -\frac{(12)^{1/3}}{2} $} & $-\seqnum{A319034}$    \\
$\bf 167$ & $  $  & $   $  &  $       $    &  $(-)$ & {\color{blue}$\bf +\frac{(12)^{1/3}}{2} $} & $+\seqnum{A319034}$    \\
{\bf tbc} &  $\bf \cdots $  & $\bf \cdots $  &  $\bf \cdots $ &  $\bf \cdots $ &  $\bf \cdots $  & $\bf \cdots $ \\
\hline
\end{tabular}
\end{center}
\vfill
\eject
\begin{center}
{\large {\bf Table 5: Continued: Real algebraic numbers $\bf \{{\boldsymbol \omega}(c)\}_{c=168}^{211}$ of height $\bf n$ and degree $\bf k$.}}
\par\smallskip\noindent
\begin{tabular}{|r|r|r|c|c|c|c|}\hline
&& && &&\\
$\bf c$ &  $\bf height\, n$ & $\bf degree\, k$ & \bf Composition & $\bf signs $  & $\boldsymbol \omega$ & \bf A-number \\
&& && && \\ \hline\hline
$\bf 168$ &   $7$ cont.   &   $3$ cont.   & $[3,1,0,1]$  &  $(+,+)$ & {\color{blue}$\bf -0.8241226211... $} & $-\seqnum{A357465}$ \\
$\bf 169$ &          &         &              &  $(-,-)$ & {\color{blue}$\bf +0.8241226211...$} & $+\seqnum{A357465}$ \\ 
$\bf 170$ &          &         &              &  $(-,+)$ & {\color{blue}$\bf -0.5981934981... $} & $-\seqnum{A357464}$ \\
$\bf 171$ &          &         &              &  $(+,-)$ & {\color{blue}$\bf +0.5981934981...$} & $+\seqnum{A357464}$ \\ 
$\bf 172$ &          &         & $ [3,0,1,1]$ &  $(+,+)$ & {\color{blue}$\bf -0.5365651646... $} & $-\seqnum{A357467}$ \\
$\bf 173$ &          &         &              &  $(+,-)$ & {\color{blue}$\bf +0.5365651646...$} & $+\seqnum{A357467}$ \\ 
$\bf 174$ &          &         &              &  $(-,+)$ & {\color{blue}$\bf -0.8513830728... $} & $-\seqnum{A357466}$ \\
$\bf 175$ &          &         &              &  $(-,-)$ & {\color{blue}$\bf +.0.8513830728..$} & $+\seqnum{A357466}$ \\ 
$\bf 176$ &      &        & $[1,3,0,1]$  &  $(+,+)$ & {\color{blue}$\bf -3.1038034027...$} & $-(\seqnum{A357103}\sspp 1)$ \\ 
$\bf 177$ &      &       & $        $    &  $(-,-) $ & {\color{blue}$\bf +3.1038034027... $}& $+(\seqnum{A357103}\sspp 1)$ \\
$\bf 178$ &      &        & $   $  &    $(+,-)$ & {\color{blue}$\bf -2.8793852415... $} & $-(\seqnum{A332437} \sspp 1)$ \\
$\bf 179$ &      &       & $        $  &  $     $ & {\color{blue}$\bf -0.6527036446... $} & $-\seqnum{A178959}$ \\
$\bf 180$ &      &       & $   $   & $    $ & {\color{blue}$\bf +0.5320888862...$} & $+(\seqnum{A322438} \sspm 3$) \\
$\bf 181$ &      &       & $   $   &   $ (-,+)$  & {\color{blue}$\bf -0.5320888862...$} & $-(\seqnum{A322438} \sspm 3$) \\
$\bf 182$ &      &       & $        $  &  $     $ & {\color{blue}$\bf +0.6527036446... $} & $+\seqnum{A178959}$ \\
$\bf 183$ &      &        & $        $  &  $  $ & {\color{blue}$\bf +2.8793852415... $} & $+(\seqnum{A332437} \sspp 1)$ \\
$\bf 184$ &      &        & $ [1,0,3,1]  $  &    $(-,+)$ & {\color{blue}$\bf -1.8793852415... $} & $-\seqnum{A332437}$ \\
$\bf 185$ &      &       & $        $  &  $     $ & {\color{blue}$\bf +0.3472963553... $} & $+\seqnum{A130880}$ \\
$\bf 186$ &      &       & $   $   & $    $ & {\color{blue}$\bf +1.5320888862... $} & $+(\seqnum{A332438} \sspm 2$) \\
$\bf 187$ &      &       & $   $   &   $ (-,-)$  & {\color{blue}$\bf -1.5320888862...$} & $-(\seqnum{A332438} \sspm 2$) \\
$\bf 188$ &      &       & $        $  &  $     $ & {\color{blue}$\bf -0.3472963553... $} & $-\seqnum{A130880}$ \\
$\bf 189$ &      &        & $        $  &  $  $ & {\color{blue}$\bf +1.8793852415... $} & $+\seqnum{A332437}$ \\
$\bf 190$ &      &        & $       $  &  $(+,+)$  & {\color{blue}$\bf -0.3221853546...$} & $-\seqnum{A357104}$ \\ 
$\bf 191$ &      &       & $        $    &  $(+,-) $ & {\color{blue}$\bf +0.3221853546... $}& $+\seqnum{A357104} $ \\
$\bf 192$ &      &        & $[1,1,0,3]$  &  $(+,+)$  & {\color{blue}$\bf -1.8637065278...$} & $-\seqnum{A356034}$ \\ 
$\bf 193$ &      &       & $        $    &  $(-,-) $ & {\color{blue}$\bf +1.8637065278... $}& $+\seqnum{A356034} $ \\
$\bf 194$ &      &        & $       $    &  $(-,+)$  & {\color{blue}$\bf -1.1745594102...$} & $-\seqnum{A357100}$ \\ 
$\bf 195$ &      &       & $        $    &  $(+,-) $ & {\color{blue}$\bf +1.1745594102... $}& $+\seqnum{A357100} $ \\
$\bf 196$ &      &        & $[1,0,1,3]$  &  $(-,+)$  & {\color{blue}$\bf -1.6716998816...$} & $-\seqnum{A294644}$ \\ 
$\bf 197$ &      &       & $        $    &  $(-,-) $ & {\color{blue}$\bf +1.6716998816... $}& $+\seqnum{A294644} $ \\
$\bf 198$ &      &        & $       $    &  $(+,+)$  & {\color{blue}$\bf -1.2134116627...$} & $-\seqnum{A337569}$ \\ 
$\bf 199$ &      &       & $        $    &  $(+,-) $ & {\color{blue}$\bf +1.2134116627... $}& $+\seqnum{A337569} $ \\
$\bf 200$ &      &        & $[2,2,0,1]$  &  $(+,+)$ & {\color{blue}$\bf -1.2971565081...$} & $-\seqnum{A357109}$ \\ 
$\bf 201$ &      &       & $        $    &  $(-,-) $ & {\color{blue}$\bf +1.2971565081... $}& $+\seqnum{A357109}$ \\
$\bf 202$ &      &        & $$           &  $(-,+)$ & {\color{blue}$\bf -0.5651977173...$} & $-\seqnum{A273065}$ \\ 
$\bf 203$ &      &       & $        $    &  $(+,-) $ & {\color{blue}$\bf +0.5651977173... $}& $+\seqnum{A273065}$ \\
$\bf 204$ &      &        & $[2,0,2,1]$  &  $(-,+)$ & {\color{blue}$\bf -1.1914878839...$} & $-(\seqnum{A316711}\sspm 1)$ \\ 
$\bf 205$ &      &       & $        $    &  $(-,-) $ & {\color{blue}$\bf +1.1914878839... $}& $+(\seqnum{A316711}\sspm 1)$ \\
$\bf 206$ &      &        & $$           &  $(+,+)$ & {\color{blue}$\bf -0.4238537990...$} & $-\seqnum{A357463}$ \\ 
$\bf 207$ &      &       & $        $    &  $(+,-) $ & {\color{blue}$\bf +0.4238537990... $}& $+\seqnum{A357463}$ \\
$\bf 208$ &      &        & $[2,1,0,2]$  &  $(+,+)$  & {\color{blue}$\bf -1.1974293369...$} & $-\seqnum{A357105}$ \\ 
$\bf 209$ &      &       & $        $    &  $(-,-) $ & {\color{blue}$\bf +1.1974293369... $}& $+\seqnum{A357105}$ \\
$\bf 210$ &      &        & $$           &  $(-,+)$  & {\color{blue}$\bf -0.8580943294...$} & $-\seqnum{A357106}$ \\ 
$\bf 211$ &      &       & $        $    &  $(+,-) $ & {\color{blue}$\bf +0.8580943294... $}& $+\seqnum{A357106}$ \\
{\bf tbc} &  $\bf \cdots $  & $\bf \cdots $  &  $\bf \cdots $   &  $\bf \cdots $ &  $\bf \cdots $  & $\bf \cdots $ \\
\hline
\end{tabular}
\end{center}
\vfill
\eject
\begin{center}
{\large {\bf Table 6: Continued: Real algebraic numbers $\bf \{{\boldsymbol \omega}(c)\}_{c=212}^{257}$ of height $\bf n$ and degree $\bf k$.}}
\par\smallskip\noindent
\begin{tabular}{|r|r|r|c|c|c|c|}\hline
&& && &&\\
$\bf c$ &  $\bf height\, n$ & $\bf degree\, k$ & \bf Composition & $\bf signs $  & $\boldsymbol \omega$ & \bf A-number \\
\hline\hline
$\bf 212$ &   $7$ cont.  &    $3$ cont.   & $[2,0,1,2]$  &  $(-,+)$  & {\color{blue}$\bf -1.1653730430...$} & $-\seqnum{A357107}$ \\ 
$\bf 213$ &      &       & $        $      &  $(-,-) $ & {\color{blue}$\bf +1.1653730430... $}& $+\seqnum{A357107}$ \\
$\bf 214$ &      &        & $$             &  $(+,+)$  & {\color{blue}$\bf -0.8351223484...$} & $-\seqnum{A357108}$ \\ 
$\bf 215$ &      &       & $        $      &  $(+,-) $ & {\color{blue}$\bf +0.8351223484... $}& $+\seqnum{A357108}$ \\
$\bf 216$ &      &       & $[1,2,0,2]$  &  $(+,+)$  & {\color{blue}$\bf -2.3593040859...$} & $-\seqnum{A357101}$ \\ 
$\bf 217$ &      &       & $        $   &  $(-,-) $ & {\color{blue}$\bf +2.3593040859... $}& $+\seqnum{A357101}$ \\
$\bf 218$ &      &        & $$          &  $(-,+)$  & {\color{blue}$\bf -0.8392867552...$} & $-(\seqnum{A058265}\sspm 1)$ \\ 
$\bf 219$ &      &       & $        $    &  $(+,-) $ & {\color{blue}$\bf +0.8392867552... $}& $+(\seqnum{A058265}\sspm 1)$ \\
$\bf 220$ &      &       & $[1,0,2,2]$   &  $(-,+)$  & {\color{blue}$\bf -1.7692923542...$} & $-\seqnum{A273066}$ \\ 
$\bf 221$ &      &       & $        $    &  $(-,-) $ & {\color{blue}$\bf +1.7692923542... $}& $+\seqnum{A273066}$ \\
$\bf 222$ &      &        & $$           &  $(+,+)$  & {\color{blue}$\bf -0.7709169970...$} & $-\seqnum{A357102}$ \\ 
$\bf 223$ &      &       & $        $    &  $(+,-) $ & {\color{blue}$\bf +0.7709169970... $}& $+\seqnum{A357102}$ \\
$\bf 224$ &      &       & $[2,1,1,1]$   &  $(+,-,+)$  & {\color{blue}$\bf -1.2337519285...$} & $-\seqnum{A358182}$ \\
$\bf 225$ &      &       & $         $   &  $(-,-,-)$  & {\color{blue}$\bf +1.2337519285...$} & $+\seqnum{A358182}$ \\
$\bf 226$ &      &       & $         $   &  $(-,-,+)$  & {\color{blue}$\bf -0.8294835409...$} & $-\seqnum{A358183}$ \\
$\bf 227$ &      &       & $         $   &  $(+,-,-)$  & {\color{blue}$\bf +0.8294835409...$} & $+\seqnum{A358183}$ \\
$\bf 228$ &      &       & $         $   &  $(+,+,+)$  & {\color{blue}$\bf -0.7389836215...$} & $-\seqnum{A358184}$ \\
$\bf 229$ &      &       & $         $   &  $(-,+,-)$  & {\color{blue}$\bf +0.7389836215...$} & $+\seqnum{A358184}$ \\
$\bf 230$ &      &       & $[1,2,1,1]$   &  $(+,-,+)$  & {\color{blue}$\bf -2.5468182768...$} & $-\seqnum{A358181}$ \\
$\bf 231$ &      &       & $         $   &  $(-,-,-)$  & {\color{blue}$\bf +2.5468182768...$} & $+\seqnum{A358181}$ \\
$\bf 232$ &      &       & $         $   &  $(+,-,-)$  & {\color{blue}$\bf -2.2469796037...$} & $-\seqnum{A231187}$ \\
$\bf 233$ &      &       & $         $   &  $       $  & {\color{blue}$\bf -0.5549581320...$} & $-\seqnum{A255240}$ \\
$\bf 234$ &      &       & $         $   &  $       $  & {\color{blue}$\bf +0.8019377358...$} & $+(\seqnum{A160389}\sspm 1)$ \\
$\bf 235$ &      &       & $         $   &  $(+,+,+)$  & {\color{blue}$\bf -1.7548776662...$} & $-\seqnum{A109134}$ \\
$\bf 236$ &      &       & $         $   &  $(-,+,-)$  & {\color{blue}$\bf +1.7548776662...$} & $+\seqnum{A109134}$ \\
$\bf 237$ &      &       & $         $   &  $(-,-,+)$  & {\color{blue}$\bf -0.8019377358...$} & $-(\seqnum{A160389}\sspm 1)$ \\
$\bf 238$ &      &       & $         $   &  $       $  & {\color{blue}$\bf +0.5549581320...$} & $+\seqnum{A255240}$ \\
$\bf 239$ &      &       & $         $   &  $       $  & {\color{blue}$\bf +2.2469796037...$} & $+\seqnum{A231187}$ \\
$\bf 240$ &      &       & $         $   &  $(-,+,+)$  & {\color{blue}$\bf -0.4655712318...$} & $-\seqnum{A088559}$ \\
$\bf 241$ &      &       & $         $   &  $(+,+,-)$  & {\color{blue}$\bf +0.4655712318...$} & $+\seqnum{A088559}$ \\
$\bf 242$ &      &       & $[1,1,2,1]$   &  $(+,-,+)$  & {\color{blue}$\bf -2.1478990357...$} & $-\seqnum{A357470}$ \\
$\bf 243$ &      &       & $         $   &  $(-,-,-)$  & {\color{blue}$\bf +2.1478990357...$} & $+\seqnum{A357470}$ \\
$\bf 244$ &      &       & $         $   &  $(+,-,-)$  & {\color{blue}$\bf -1.8019377358...$} & $-\seqnum{A160389}$ \\
$\bf 245$ &      &       & $         $   &  $       $  & {\color{blue}$\bf -0.4450418679...$} & $-\seqnum{A255241}$ \\
$\bf 246$ &      &       & $         $   &  $       $  & {\color{blue}$\bf +1.2469796037...$} & $+\seqnum{A225249}$ \\
$\bf 247$ &      &       & $         $   &  $(-,-,+)$  & {\color{blue}$\bf -1.2469796037...$} & $-\seqnum{A225249}$ \\
$\bf 248$ &      &       & $         $   &  $       $  & {\color{blue}$\bf +0.4450418679...$} & $+\seqnum{A255241}$ \\
$\bf 249$ &      &       & $         $   &  $       $  & {\color{blue}$\bf +1.8019377358...$} & $+\seqnum{A160389}$ \\
$\bf 250$ &      &       & $         $   &  $(+,+,+)$  & {\color{blue}$\bf -0.5698402909...$} & $-\seqnum{A357471}$ \\
$\bf 251$ &      &       & $         $   &  $(-,+,-)$  & {\color{blue}$\bf +0.5698402909...$} & $+\seqnum{A357471}$ \\
$\bf 252$ &      &       & $         $   &  $(-,+,+)$  & {\color{blue}$\bf -0.3926467817...$} & $-\seqnum{A357472}$ \\
$\bf 253$ &      &       & $         $   &  $(+,+,-)$  & {\color{blue}$\bf +0.3926467817...$} & $+\seqnum{A357472}$ \\
$\bf 254$ &      &       & $ [1,1,1,2]$  &  $(+,+,+)$  & {\color{blue}$\bf -1.3532099641...$} & $-\seqnum{A357469}$ \\
$\bf 255$ &      &       & $         $   &  $(-,+,-)$  & {\color{blue}$\bf +1.3532099641...$} & $+\seqnum{A357469}$ \\
$\bf 256$ &      &       & $         $  &  $(-,-,+)$  & {\color{blue}$\bf -1.2055694304...$} & $-\seqnum{A137421}$ \\
$\bf 257$ &      &       & $         $   &  $(+,-,-)$  & {\color{blue}$\bf +1.2055694304...$} & $+\seqnum{A137421}$ \\
{\bf tbc} &  $\bf \cdots $  & $\bf \cdots $  &  $\bf \cdots $   &  $\bf \cdots $ &  $\bf \cdots $  & $\bf \cdots $ \\
\hline
\end{tabular}
\end{center}
\vfill
\eject
\begin{center}
{\large {\bf Table 7: Continued: Real algebraic numbers $\bf {\{\boldsymbol\omega}(c)\}_{258}^{291}$ of height $\bf n$ and degree $\bf k$.}}
\par\smallskip\noindent
\begin{tabular}{|r|r|r|c|c|c|c|}\hline
&& && &&\\
$\bf c$ &  $\bf height\, n$ & $\bf degree\, k$ & \bf Composition & $\bf signs $  & $\boldsymbol \omega$  & \bf A-number \\
&& && && \\ \hline\hline
$\bf 258$ &   $7$ cont.   &  $3$ cont.    & $ [1,1,1,2]$ cont.  &  $(-,+,+)$  & {\color{blue}$\bf -0.8105357137...$} & $-\seqnum{A357468}$ \\
$\bf 259$ &         &          & $                $   &  $(+,+,-)$  & {\color{blue}$\bf +0.8105357137...$} & $+\seqnum{A357468}$ \\
\hline
$\bf 260$ &   $  $  & $ 4 $ & $[3,1]$  &  $(-)$ & {\color{blue}$\bf -\frac{27^{\frac{1}{4}}}{3}$} & $-\seqnum{A358186}$    \\
$\bf 261$ &   $  $  & $   $ & $    $   &        & {\color{blue}$\bf  +\frac{27^{\frac{1}{4}}}{3}$} & $+\seqnum{A358186}$    \\
$\bf 262$ &   $  $  & $  $   & $[1,3]$  &  $(-)$ & {\color{blue}$\bf -3^{\frac{1}
{4}}$} & $-\seqnum{A011002}$    \\
$\bf 263$ &   $  $  & $  $   & $     $  &        & {\color{blue}$\bf +3^{\frac{1}
{4}}$} & $+\seqnum{A011002}$    \\
$\bf 264$ &   $  $  & $  $ & $[1,2,0,0,1]$  &  $(+,-)$ & {\color{blue}$\bf -2.1069193403...$} & $-\seqnum{A358188}$    \\
$\bf 265$ &   $  $  & $  $ & $           $  &  $    $ & {\color{blue} $\bf +0.7166727492...$} & $+\seqnum{A358187}$    \\
$\bf 266$ &   $  $  & $  $ & $           $  &  $(-,-)$ & {\color{blue} $\bf -0.7166727492...$} & $-\seqnum{A358187}$    \\
$\bf 267$ &   $  $  & $  $ & $           $  &  $     $ & {\color{blue}$\bf +2.1069193403...$} & $+\seqnum{A358188}$    \\
$\bf 268$ &   $  $  & $  $ & $[1,0,2,0,1]$  &  $(-,-)$ & {\color{blue}$\bf -1.5537739740...$} & $-\seqnum{A278928}$ \\
$\bf 269$ &   $  $  & $   $ & $           $  &  $     $ & {\color{blue}$\bf +1.5537739740...$} & $+\seqnum{A278928}$ \\
$\bf 270$ &   $  $  & $   $ & $           $  &  $(+,-)$ & {\color{blue}$\bf -0.6435942529...$} & $-\seqnum{A154747}$ \\
$\bf 271$ &   $  $  & $   $ & $           $  &  $     $ & {\color{blue}$\bf +0.6435942529... $} & $+\seqnum{A154747}$  \\
$\bf 272$ &   $  $  & $   $ & $[1,0,0,2,1]$  &  $(+,-)$ & {\color{blue}$\bf -1.3953369944...$} & $-\seqnum{A358190}$ \\
$\bf 273$ &   $  $  & $   $ & $           $  &  $     $ & {\color{blue}$\bf +0.4746266175...$} & $+\seqnum{A358189}$ \\
$\bf 274$ &   $  $  & $   $ & $           $  &  $(-,-)$ & {\color{blue}$\bf -0.4746266175...$} & $-\seqnum{A358189 }$ \\
$\bf 275$ &   $  $  & $   $ & $           $  &  $     $ & {\color{blue}$\bf +1.3953369944... $} & $+\seqnum{A358190}$  \\
\hline 
$\bf 276$ &   $  $  & $ 5 $ & $[2,1]$  &  $(+)$ & {\color{blue}$\bf -{\frac{1}{2}}^{\frac{1}{5}} $} & $-\seqnum{A358938}$    \\
$\bf 277$ &   $  $  & $   $ & $     $  &  $(-)$ & {\color{blue}$\bf +{\frac{1}{2}}^{\frac{1}{5}} $} & $+\seqnum{A358938}$    \\
$\bf 278$ &   $  $  & $   $ & $[1,2]$  &  $(+)$ & {\color{blue}$\bf -2^{\frac{1}{5}}$} & $-\seqnum{A005531}$    \\
$\bf 279$ &   $  $  & $   $ & $    $  &  $(-)$ & {\color{blue}$\bf  +2^{\frac{1}{5}}$} & $+\seqnum{A005531}$    \\
$\bf 280$ &   $  $  & $   $ & $[1,1,0,0,0,1]$  &  $(-,+)$ & {\color{blue}$\bf -0.8566748838... $} & $ -\seqnum{A230152}$    \\
$\bf 281$ &   $  $  & $   $ & $             $  &  $(+,-)$ & {\color{blue}$\bf +0.8566748838... $} & $ +\seqnum{A230152}$   \\
$\bf 282$ &   $  $  & $   $ & $[1,0,1,0,0,1]$  &  $(-,+)$ & {\color{blue}$\bf -1.2365057033... $} & $-\seqnum{A358940}$    \\
$\bf 283$ &   $  $  & $   $ & $             $  &  $(-,-)$ & {\color{blue}$\bf +1.2365057033... $} & $ +\seqnum{A358940}$   \\
$\bf 284$ &   $  $  & $   $ & $             $  &  $(+,+)$ & {\color{blue}$\bf -0.8376197748... $} & $-\seqnum{A358939}$    \\
$\bf 285$ &   $  $  & $   $ & $             $  &  $(+,-)$ & {\color{blue}$\bf +0.8376197748... $} & $ +\seqnum{A358939}$   \\
$\bf 286$ &   $  $  & $   $ & $[1,0,0,1,0,1]$  &  $(+,+)$ & {\color{blue}$\bf -1.1938591113... $} & $-\seqnum{A358942}$    \\
$\bf 287$ &   $  $  & $   $ & $             $  &  $(-,-)$ & {\color{blue}$\bf +1.1938591113... $} & $ +\seqnum{A358942}$   \\
$\bf 288$ &   $  $  & $   $ & $             $  &  $(-,+)$ & {\color{blue}$\bf -0.8087306004... $} & $-\seqnum{A358941}$    \\
$\bf 289$ &   $  $  & $   $ & $             $  &  $(+,-)$ & {\color{blue}$\bf +0.8087306004... $} & $ +\seqnum{A358941}$   \\
$\bf 290$ &   $  $  & $   $ & $[1,0,0,0,1,1]$  &  $(-,+)$ & {\color{blue}$\bf -1.1673039782... $} & $-\seqnum{A160155}$    \\
$\bf 291$ &   $  $  & $   $ & $             $  &  $(-,-)$ & {\color{blue}$\bf +1.1673039782... $} & $ +\seqnum{A160155}$   \\
\hline
$\bf \cdots $ &  $\bf \cdots $  & $\bf \cdots $  &  $\bf \cdots $   &  $\bf \cdots $ &  $\bf \cdots $  & $\bf \cdots $ \\
\hline
\end{tabular}
\end{center}
\vfill
\eject
\begin{center}
{\large {\bf Table 9: $\bf Phi(n,\,k)$, the number of real algebraic numbers $\boldsymbol\omega$ of height $\bf n$ and degree $\bf k$.}}
\par\smallskip\noindent
\begin{tabular}{|l|c|c|c|c|c|c|c|}\hline
&& && && &\\
$\bf n\, \setminus \, k$ &  $\bf 1$ & $\bf2$ & $\bf3$& $\bf4$ & $\bf5$ & $\bf6\, ...$ & $\bf Phi(n)$ \\
&& && && & \\ \hline\hline
$\bf 1$ &   $ 1 $  &  & & & & & $ 1 $  \\
\hline
$\bf 2$ &    $ 2 $  &  & & & & & $2 $\\
\hline
$\bf 3$ &    $ 4  $  & $ 0 $ &   & & & & $4$  \\
\hline
$\bf 4$ &    $ 4 $  &  $ 8 $ & $ 0$ & & & & $12$    \\
\hline
$\bf 5$ &    $ 8 $  &  $ 8 $ & $12$ & $ 0 $ & & & $28$  \\
\hline
$\bf 6$ &    $ 4 $  & $ 32 $ & $ 20 $ &  $16$ & $0$ &  & $72$   \\
\hline
$\bf 7$ &   $ 12 $  &  $ 28 $ & $ 100 $&  $16 $ & $ 16$ & $0$ & $172$  \\
\hline
$\bf \cdots $  &  $\bf \cdots $   &  $\bf \cdots $ &  $\bf \cdots $  & $\bf \cdots $ & $\bf \cdots $ & $\bf \cdots $ &  $\bf \cdots$ \\
\hline
\hline
\end{tabular}
\end{center}
\vfill
\eject

\end{document}